\documentclass[12pt]{article}
\usepackage{amsmath,amssymb,latexsym,graphicx}
\newtheorem{theorem}{Theorem}
\newtheorem{lemma}[theorem]{Lemma}
\newtheorem{proposition}[theorem]{Proposition}
\newtheorem{remark}[theorem]{Remark}
\newtheorem{corollary}[theorem]{Corollary}
\newtheorem{definition}[theorem]{Definition}
\newtheorem{conjecture}[theorem]{Conjecture}

\def\pf{{\bf Proof }}

\begin{document}
\title{Boundary Forelli theorem for the sphere in $\mathbb C^n$ and $n+1$ bundles of complex lines}
\author{Mark~ L.~Agranovsky}
\maketitle
\begin{abstract}
Let $B^n$ be the unit ball in $\mathbb C^n$ and let the points $a_1,\cdots,a_{n+1} \in  B^n$ are affinely independent. 
If $f \in C(\partial B^n)$ and for any complex line
$L,$ containg at least one of the point $a_j$ the restriction $f\vert_{L \cap \partial B^n}$
extends holomorphically in the disc $L \cap B^n$, then $f$ is the boundary value of a holomorphic function in $B^n$. 
The condition for the points $a_j$ is sharp.

\end{abstract}

\section{Introduction}\label{S:Intro}

Forelli theorem \cite{R} on holomorphicity on complex slices, which, in a sense, can be regarded as a variation of classical Hartogs' theorem about separate analyicity,
 says that if a function $f$ of $n$ complex variables in a
domain $D \subset \mathbb C^n$ is holomorphic on each complex line $L$ passing through a fixed point $a \in D$  and $f \in C^{\infty}(\{a\})$ then $f$ is holomorphic in a neighborhood of $a$ (in fact, in a maximal subdomain which is complete circular with respect to $a$).

One might think of a direct boundary analog of Forelli theorem obtained by replacing the words "`holomorphic function"'
by "boundary value of holomorphic function"' and "`holomorphic on a complex line $L$ "' by "`holomorphically extends in $L$"'.

However, such direct analog fails. Indeed, let $B^n$ be the unit ball in $\mathbb C^n.$ The function 
$f(z)=|z_1|^2$ on the unit sphere $\partial B^n=\{|z_1^2+\cdots |z_n|^2=1\}$  is constant on circles  $L \cap \partial B^n$ obtained as intersections of the unit sphere with the complex lines $L$
passing through $0.$ Thus, $f$ extends holomorphically (as a constant) in the discs $L \cap B^n$ but does not extend in $B^n$ as a holomorphic function of $n$ variables,  because it is real-valued and non-constant.

In the recent preprint \cite{APrep}, the author proved that if $f$ is real analytic on the complex sphere
$\partial B^n$ then two bundles of complex lines suffice. More precisely, the result of \cite{APrep} says that
if $f \in C^{\omega}(\partial B^n)$ and  $a,b \in B^n, \ a \neq b,$ 
and if for each complex line $L,$ containing at least one of the points $a,b,$ the 
function $f$ holomorphically extends in each disc $L \cap \partial B^n,$ then $f$ extends in $B^n$ as a holomorphic function.
J. Globevnik \cite{Gl} strenghened this result, replacing $C^{\omega}$ by $C^{\infty}.$ He also studied in details the case when one of the point $a,b$ or both of them belong to the 
interior of the ball.  L. Baracco \cite{B} gave an alternative proof of the result in \cite{APrep}, which extends also to so called extremal discs.

Regularity of $f$ is essential. The following example was found by J. Globevnik. 
The function $f(z)=z_n^k /\overline z_n$ on the unit sphere $\partial B^n$
extends holomorphically through each complex line intersecting the 
hyperplane $z_n=0$ but does not admit holomorphic extension in the ball $B^n.$
Notice for any prescribed $r$  the function $f$ belongs to $ C^r(\partial B^n)$ for sufficiently large  $k$, 
but $f \notin C^{\infty}(\partial B^n).$   

In the above example, the vertices of the bundles of the complex lines lie on a hyperplane. 
It was conjectured in \cite {APrep} that this is the only obstruction for boundary Forelli theorem to be true for continuous functions: 
\begin{conjecture} (\cite{APrep})Any $n+1$ bundles of complex lines with the vertices $a_1,...,a_{n+1} \in B^n$ in general position do test boundary values $f \in C(\partial B^n)$ of holomorphic functions
in the unit ball $B^n \subset \mathbb C^n.$ 
\end{conjecture}

In this article, we prove this conjecture (Theorem \ref{T:main}). 
The proof is based on a modification of ideas and constructions from \cite{APrep}, where the reduction was used to a problem of characterization of polyanalytic functions in the unit disc \cite{Agr}. 

For some earlier results in the problem, we  refer to articles \cite{AS},\cite{GS},\cite{St}. More extended bibliography can be found, for instance, in \cite{APrep}.   
\section{Main results}

Let us formulate our main result.
\begin{theorem}\label{T:main}
Let $A=\{a_1,\cdots,a_{n+1}\} \subset  B^n$  be $n+1$ points in general position, i.e. belonging to no $(n-1)$-dimensional complex plane. Let $f \in C(\partial B^n)$
and assume that the following condition holds

(*) for every complex line $L \subset \mathbb C^n$ such that $L \cap A \neq \emptyset$ there exists a function
$F_L \in C(L \cap \overline B^n),$ holomorphic in $L \cap B^n$ and such that $F_L(z)=f(z)$ for $z \in L \cap \partial B^n.$ 

Then there exists a function $F \in C(\overline B^n),$ holomorphic in $B^n$ ans such that $F(z)=f(z)$ for 
$z \in \partial B^n.$ 
\end{theorem}
\begin{remark} Obviously, the set $A$ in Theorem \ref{T:main} can be taken any, not necessary finite,nonempty subset of $B^n$  belonging to no complex hyperplane. 
\end{remark}

Denote $A(\partial B^n)$ the ball-algebra of all continuous functions that extend continuously in $\overline B^n$ as a function holomorophic in the open ball.

Given $a \in \mathbb C^n$, denote  
$$\mathcal L_a= \{ \mbox{all complex lines in} \ \mathbb C^n \ \mbox{containing} \  a \}.$$

In this notations,  Theorem \ref {T:main} asserts that
the holomorphic extendibility of $f$ along complex lines $L \in \mathcal L_{a_1} \cup \cdots \cup \mathcal L_{a_{n+1}},$ where the points $a_j$ are affinely independent, implies that $f \in A(\partial B^n).$

The proof of Theorem \label{main} rests on the following
\begin{theorem}\label{T:mainmain} Let $a,b \in B^n, \ a \neq b.$ If a function $f \in C(\partial B^n)$ extends holomorphically in every complex line
$L \in \mathcal L_a \cup \mathcal L_b$ (*) then it extends in complex lines $L \in \mathcal L_c$ for any point $c \in B^n$ on the complex line $L_{a,b}$ joining $a$ and $b$.
\end{theorem}
\begin{proposition} \label{P:reduction}
Theorem \ref{T:main} follows from Theorem \ref{T:mainmain} for $n=2.$
\end{proposition}
\pf
Theorem \ref{T:mainmain} implies that $f$ extends from complex lines meeting the affine span of the set $A=\{a_1,\cdots,a_{n+1}\}.$
By the condition, $span A=\mathbb C^n$ and therefore $f$ extends in all cross-sections $L \cap \partial B^n$ by arbitrary complex line $L.$ 
Then it is easy to show that $f \in A(\partial B^n)$.
For instance, by convolving on the unitary group $SU(2)$, one can approximate $f$ by smooth functions. Applying Cauchy-Green formula to cross-sections $L \cap B^n$ and
letting $L$ to tend to a tangent line at the sphere $\partial B^n$ lead to tangential CR-conditions for $f$ on the boundary. This implies $f \in A(\partial B^n).$

\begin{proposition}\label{P:n=2}
The case of arbitary $n$ in Theorem \ref{T:mainmain} follows from the case $n=2.$ 
\end{proposition}
\pf
Take an arbitrary point $c \in L_{a,b} \cap B^n$ and an arbitary complex line $L_c$ containing $c.$
Let $\Pi$ be the complex 2-plane spanned by the complex lines $L_{a,b}$ and $L_c.$ Then $\Pi \cap B^n$ is a complex ball and the restriction $f \vert_{\Pi \cap \partial B^n}$ extends in each complex
line $L \subset \Pi$ passing through $a$ or $b.$ Therefore, we are in the situation of Theorem \ref{T:mainmain} and by assumption $f $ extends holomorphically in the complex line $L_c.$
Thus, $f$ extends holomorphically in complex lines from $\mathcal L_c$ and Theorem \ref{T:mainmain} in arbitrary dimension follows.

In fact, the proof of Theorem \ref{T:mainmain} contains an explicit description of functions $f \in C(\partial B^n)$ that extend in the complex lines from two bundles $\mathcal L_a$ and $\mathcal L_b$, $a,b \in B^2.$
For simplicity, we restrict ourselve by the case $n=2.$

Applying suitable automorphism form $Aut(B^2)$ we can assume, without loss of generality that $a=(a_1,0), b=(b_1,0).$
\begin{theorem}\label{T:mainmainmain} Let $a, b \in B^2 \cap \{z_2=0\}.$ and $f \in C(\partial B^2).$ Then $f$ extends holomorphically from the complex lines, passing through at least one of the  points $a$ and $b$, 
if and only if Fourier series in the angular variable of $z_2$ has the form
\begin{equation}\label{E:fur}
f(z_1,re^{i\phi})=\sum\limits_{\nu=0}^{\infty}F^{\nu}(z_1,r)r^{\nu}e^{i\nu\phi}=\sum\limits_{\nu}^{\infty}F^{\nu}(z_1,|z_2|)z_2^{\nu}, \ z_2=re^{i\phi},
\end{equation}
and  the coefficients $F^{\nu}$ are representable as follows:
\begin{equation}\label{E:ff}
F^{\nu}(z_1,|z_2|)=\sum\limits_{0 \leq 2j<\nu}\frac{h^{\nu}_{j}(z_1)}{|z_2|^{2j}},
\end{equation}
where $h_j$ are continuous functions in the disc $\{|z_1| \leq 1\},$ holomorphic in $\{|z_1|<1\}.$
\end{theorem}
Clearly, the above mentioned example of J. Globevnik is a particular case of the functions described by Theorem \ref{T:mainmainmain}. 
\begin{corollary}(\cite{APrep},\cite{Gl}) If $f \in C^{\infty}(\partial B^n)$ then extendibility in the complex lines only from two bundles $\mathcal L_a, \mathcal L_b,$ where $a,b \in B^n, \ a\neq b$, implies
$f \in A(\partial B^n).$
\end{corollary}
\pf
First of all, it is shown in \cite{APrep} that it suffices to prove the assertion for the case $n=2.$

The condition $f \in C^{\infty}(\partial B^2)$ demands that all the terms in the representation (\ref{E:fur}),(\ref{E:ff}):
$$f(z_1,z_2)=\sum\limits_{\nu=0}^{\infty} \sum\limits_{0 \leq 2j<\nu}h^{\nu}_{j}(z_1)\frac{z_2^\nu}{|z_2|^{2j}},$$
with $j>0$ 
must be be identically zero. 
In other words,
$$F^{\nu}(z_1)=h^{\nu}_0(z_1)$$ 
is holomorphic. 
Then
\begin{equation}\label{E:holo}
f(z_1,z_2)=\sum\limits_{\nu=0}^{\infty}h^{\nu}_0(z_1)z_2^{\nu}, \ |z_1|^2+|z_2|^2=1,
\end{equation}
extends holomorphically in $|z_1|^2+|z_2|^2 <1.$ The extension is delivered just by the right hand side in (\ref{E:holo}).  

\section{The approach} 

We follow the strategy of our previous articles \cite{APrep} and \cite{Agr}. There the test of holomorphic extendibilty from the complex sphere was obtained for
regular functions $f$ and two bundles of complex lines. The main case was $n=2.$ We assumed, without loss of generality, 
 that the vertices of the bundles belong to one of the complex axis. 
 
The first step in \cite{APrep} was to develop $f$ in Fourier series with respect to the rotation group around the axis and 
reduce the original problem to a problem in one complex variable. Namely, the original problem about holomorphic extenstions in two bundles of complex lines 
converts to characterization those functions in the unit disc that admit meromorphic extension into two 
families of concentric hyperbolic circles. Such problem was solved in \cite{Agr}, even for general one-parameter families of circles. However, the regularity (real-analyticty) of functions was used essentially.
In the present article, we solve the above one-dimensional problem for much more special family, than in \cite{Agr}, namely for concentric hyperbolic circles , 
however under regularity conditions which are met in our case. 

It is worth noticing that the reduction to one-dimensional case is different from that in \cite{APrep}. 
Namely, we consider the Poisson interal associated to the Laplace-Beltrami operator, invariant
with respect to the group $Aut(B^2)$ of biholomorhic automorphisms of the ball. We tranlsate the original problem in terms of invariant Poisson integral and then reduce it to a problem in the complex disc
by double expansion: first in Fourier series with respect to the rotation group about the axis
containing the vertices of the bundles, and then developing the result in Taylor series. 

The crucial advantage of exploiting the invariant Poisson integral is that it is real-analytic and, moreover, develops into power series in  $B^2.$   
This allows to solve corresponding one-dimensional problem, by reducing the order of singularity of the meromorphic extension and then using the result from \cite{Gl} for continuous functions and hyperbolic cirlces.  
The original problem for the complex sphere follows immediately from the reduced problem in the disc. 

Moreover, the constructions in the proof of Theorem \ref{T:mainmain} (which implies \ref{T:main}) lead to an explicit description (Theorem \ref{T:mainmainmain}) of 
those continuous functions on the complex sphere that admit holomorphic extension from two bundles of complex lines. 
In particular, this contains the results of \cite{APrep},\cite{Gl} about sufficiency of 
two bundles for testing smooth boundary values of holomorphic functions on the complex sphere. 

\section{Preliminaries}  
From now on, we will focus on proving Theorem \ref{T:mainmain} for $n=2$.

Denote $G=Aut(B^2)$ the group of biholomorphic automorhisms of the complex ball $B^2$. Every automorphism is a Moebius transformation and maps affine subspaces to affine subspaces.
Given a set $A \subset B^n$ denote:
$$G_{a,b}=\{\omega \in G=Aut(B^2): \omega(0) \in \{a,b\}\}.$$
For arbitrary function $g \in C(\partial B^2)$ introduce the functions (complex moments):

\begin{equation}\label{L:Fourier}
g_{m}(z):=\int\limits_{|\lambda|=1}g(\lambda z)\lambda^m d\lambda.
\end{equation}

\begin{lemma}\label{L:moments}
The condition (*) of Theorem \ref{T:mainmain} is equivalent to
\begin{equation}\label{E:moments}
(f\circ \omega)_{m}(z)=0, \ m=0,1,\cdots, 
\mbox{for all} \ z \in \partial B^2 \ \mbox{and} \mbox{for all} \  \omega \in G_{a,b}
\end{equation}
\end{lemma}
\pf 
First, consider the bundle $\mathcal L_0$ of all complex lines pasing through the origin. Then the sections $L \cap \partial B^2$ are circles 
$$C_z=\{\lambda z: \lambda \in \mathbb C\}, \ z \in \partial B^2.$$
A function $g$ extends from the cirlce $C_z$ inside the disc $|\lambda|<1$ if and only if the complex moments $g_m(z)=0, \ m=0,1,\cdots.$ 

Now, the automorphisms $\omega \in G_A$ transform
the bundle $\mathcal L_0$ in the bundle $\mathcal L_a$ or in the bundle $\mathcal L_b.$ T
he holomorphic extendibility of $f$ from the complex lines from $\mathcal L_{a}, \mathcal L_b$ is equivalent
to holomorphic extendibilty of $f \circ \omega, \ \omega \in G_{a,b}$ from the complex lines from $\mathcal L_0$ and hence Lemma \ref{L:moments} follows.
\begin{remark} In fact, it suffices that  condition (*) of Theorem \ref{T:mainmain} or equivalent condition (\ref{E:moments}) 
hold for some two automorphisms $\omega_a,\omega_b \in Aut(B^2),$ such that $\omega_a(0)=a, \omega_b(0)=b,$
since other automorophisms from $G_{a,b}$ are obtained from $\omega_a,\omega_b$ by composing with automorphisms leaving the vertex $a$ or, 
correspondigly, the vertex $b$ fixed. The latter automorphisms map the corresponding bundles $\mathcal L_a,\mathcal L_b$ onto themselves.  
\end{remark}

\section{Invariant Poisson integral}
Introduce the invariant Poisson kernel (\cite{R}): 
\begin{equation}\label{E:Poisson}
P(z,\xi)=\frac{(1-|z|^2)^2}{|1-\langle z,\xi\rangle|^4},
\end{equation}
where $\langle z,\xi \rangle= z_1 \overline \xi_1 + z_2 \overline \xi_2$ is the inner product in $\mathbb C^2.$
The invariant Poisson kernel  $P(z,\xi)$ is defined on $(\overline B^n \times \overline B^n) \setminus \{z=\xi\}.$ 
It is well related to the Moebius group and the invariance property is expressed by the identity
$$P(\omega z,\omega \xi)=P(z,\xi)$$
for all $\omega \in G=Aut(B^2).$

Given a function $f \in C(\partial B^2)$ denote
$$P(f)(z)=\int_{\partial B^2}P(z,\xi)f(\xi)dA(\xi)$$
the invariant Poisson integral. Here $dA$ is the area measure on the unit complex sphere.

The function $P(f)$ is continuous in the closed ball and has $f$ its boundary values:
$$P(f)(z)=f(z), \ z \in \partial B^2.$$  
In the open ball $B^2$, the function $F=P(f)$ is $M$-harmonic which means
$$\widetilde \Delta F=0,$$
where $\widetilde\Delta$ is the invariant Laplacian:
$$\widetilde\Delta=4(1-|z|^2)\sum\limits_{i,k=1}^2(\delta_{ik}-z_i\overline z_k)\frac{\partial^2f}{\partial z_i \partial \overline z_k}.$$
Since the operator $\widetilde \Delta$ is elliptic, the function $F$ is real-analytic in the open ball $B^2.$ This circumstance will be important in the sequel.

\begin{lemma}\label{L:condPoisson}
The condition (*) of Theorem \ref{T:main} is equivalent to 
\begin{equation}\label{E:Poisson_moments}
(Pf \circ \omega)_{m}(z):=\int\limits_{|\lambda|=1}(Pf \circ \omega)(\lambda z)\lambda^m d\lambda=0, \ m=0,1,\cdots, \ z \in B^n,
\mbox{for all} \  \omega \in G_A.
\end{equation}
\end{lemma}
\pf First of all, due to the $Aut(B^2)-$ invariance, $Pf \circ \omega=P(f \circ \omega)$ and hence $(Pf \circ \omega)_{m}=P(f\circ\omega)_{m}=P((f\circ \omega)_m).$
Clearly $P((f \circ\omega)_{m})=0$ in $B^2$ is equivalent to $(f \circ\omega)_{m}=0$ on $\partial B^2$.


\section{Separation of variables. Moment conditions in terms of Fourier coefficients $F^{\nu}$}

As above, denote $F=Pf$. 
We fix the points $a,b \in B^2$ and denote $L_{a,b}$ the complex line containing both points $a$ and $b$. Applying a unitary rotation of the ball, we can assume, without loss of generality ,
that 
$$a=(a_1,0), b=(b_1,0).$$

Our main assumption is that all the moments vanish:
\begin{equation}\label{E:momentB^2}
(F\circ \omega)_m(z)=0, \ z \in B^2, \ m=0,1,\cdots,
\end{equation}
where $\omega \in G_{a,b} \subset Aut(B^2),$ i.e. $\omega(0) \in \{a,b\}.$

Our aim is to prove that the same moment condition holds if $\omega(0) \in L_{a,b}.$
We start with simplyfying the form of the function $F(z_1,z_2)=Pf(z_1,z_2).$

Substituting in the Poisson integral the power decomposition of the Poisson kernel (\ref{E:Poisson})
$$P(z,\xi)=(1-|z|)^2\sum\limits_{n,m=0}^{\infty}(n+1)(m+1){\langle z,\xi \rangle}^n{ \overline{ \langle z,\xi \rangle}}^m,$$
and integrating in $\xi$ yields the decomposition of the Poisson integral
$$F(z)=Pf(z)=\sum\limits_{n,m,k,l=0}^{\infty} a_{n,m,k,l}z_1^nz_2^m\overline z_1^k\overline z_2^l.$$
which uniformly converges on compacts in the open ball $B^2.$ 

Fix $\nu \in \mathbb Z.$ The terms with $m-l=\nu$ have the form $a_{n,l+\nu,k,l}z_1^n \overline z_1^k |z_2|^{2l}z_2^{\nu}$ and 
 we can  rewrite the series as follows:
\begin{equation}\label{E:rewrite}
F(z_1,z_2)=\sum\limits_{\nu=-\infty}^{\infty}F^{\nu}(z_1,|z_2|) z_2^{\nu},
\end{equation}
where  
\begin{equation}\label{E:Fnudecomp}
F^{\nu}(z_1,z_2)=F^{\nu}(z_1,|z_2|):=\sum\limits_{l=0}^{\infty}A^{\nu}_l(z_1)|z_2|^{2l} 
\end{equation}
and
\begin{equation}\label{E:Fnudecompdecomp}
A^{\nu}_l(z_1)=\sum\limits_{n,k=0}^{\infty}a_{n,l+\nu, k,l}^{\nu}z_1^n \overline z_1^k.
\end{equation}
This power series converges in $\overline B^2$, uniformly on  compacts in the open ball $B^2,$ and represents the function $F^{\nu}$ which is continuous in $B^2.$
It is easy to see that $F^{\nu}(z_1,|z_2|)z_2^{\nu}$ is obtained by integration and expresses through Fourier coefficients on the circles $|z_2|=const$: 
\begin{equation}\label{E:FnuFourier}
F^{\nu}(z_1,|z_2|)z_2^{\nu}= \frac{1}{2\pi}\int\limits_{0}^{2\pi}F(z_1,e^{i\phi}z_2)e^{-i\nu\phi}d\phi.
\end{equation}
\begin{lemma}\label{L:smooth}
The function $F^{\nu}(z_1,|z_2|)(1-|z_1|^2)^{\nu}$ is continuous in the closed ball $\overline B^2.$ 
\end{lemma}
\pf
The representation  (\ref{E:FnuFourier}) 
yields that the function
$$M(z_1,z_2):=F^{\nu}(z_1,|z_2|)z_2^{\nu}$$ 
belongs to $C(\overline B^2).$ We have on the sphere $|z_1^2+|z_2|^2=1:$
$$F^{\nu}(z_1,|z_2|)(1-|z_1|^2)^{\nu}=|z_2|^{2\nu}z_2^{-\nu}M(z_1,z_2)=\overline {z_2^{\nu}}M(z_1,z_2).$$
and the right hand side is continuous on the sphere $|z_1|^2+|z_2|^2=1.$ The function $F_{\nu}(z_1,z_2)(1-|z_1|^2)^{\nu}$
is the sum of the power series (\ref{E:Fnudecomp}) and hence it is continuous in the closed ball by Abel's theorem.

Let $c \in L_{a,b}=\{z_2=0\}, \ c=(c_1,0).$ The automorphism $\omega$ sending $0$ to $c$ can be taken
$$\omega(z)=\omega_c(z)=(\frac{z_1+c_1}{1+\overline c_1 z_1},\frac{\sqrt{(1-|c_1|^2)} z_2}{1+\overline c z_1}).$$ 

\begin{lemma} \label{L:TFAE}
The following conditions for the functions $F^{\nu}$ are equivalent to the moment condition (\ref{E:momentB^2}) for the function $F$ and the automorphism $\omega=\omega_c$:
\begin{enumerate}
\item for every $\nu >0$ and for every $z \in B^2$ the function $\lambda \mapsto (F^{\nu}\circ\omega_c)(\lambda z)$ continuously extends in the disc $|\lambda|<1$ as a meromorphic function
with the single singular point $\lambda=0$ which is the pole of order at most $\nu$,
\item for every $\nu \leq 0$ and for every $z \in B^2$ the function $\lambda \mapsto (F^{\nu}\circ\omega_c)(\lambda z)$ continuosly extends in the disc $|\lambda|<1$ as a holomorphic function.
\end{enumerate}
\end{lemma}
\pf
Let us start with the case of the identical automorphism, $\omega(z)=z.$ 
We have 
$$\int\limits_{|\lambda|=1}F(\lambda z)\lambda^m d\lambda=\sum\limits_{\nu=-\infty}^{\infty}\int\limits_{|\lambda|=1} F^{\nu}(\lambda z) \lambda^{\nu} z_2^{\nu}\lambda^m d\lambda.$$
The left hand side is 0 for all $z \in B^2$ is equivalent to the each term in the right hand side to be 0, because one can arbitrarily rotate the variable $z_2$ in the right hand side.
But vanishing complex moments
$$\int\limits_{|\lambda|=1} F^{\nu}(\lambda z) \lambda^{\nu+m} d\lambda=0$$
exactly means that the function
$$\lambda \mapsto F^{\nu}(\lambda z)\lambda^{\nu}$$
extends holomorphically in the disc $|\lambda|=1$. This proves Lemma \ref{L:TFAE} for the case $\omega=id.$
If we replace $\omega=id$ by $\omega=\omega_c$ then we will have from (\ref{E:rewrite}:
$$(F\circ \omega_c)(z)=\sum\limits_{\nu=-\infty}^{\infty}(F^{\nu}\circ\omega_c)(z) \frac{(1-|c|^2)^{\nu/2} z_2^{\nu}}{(1+\overline cz_1)^{\nu}}.$$
Therefore
$$(F\circ\omega_c)^{\nu}(z)=(F^{\nu}\circ\omega_c)(z)\frac{(1-|c|^2)^{\nu/2}}{(1+c\overline z_1)^{\nu}}.$$
The statement of Lemma \ref{L:TFAE} is checked already for $(F\circ\omega)^{\nu}.$
Since $|c|<1$, the factor $(1+\overline z_1)^{-\nu}$ is holomorphic in $|z_1| \leq 1$ and has no zeros, therefore $F^{\nu}\circ\omega_c$ possesses same type of meromorphic or analytic extendibilty.

\section{Moment conditions in terms of Taylor coefficients $A_l^{\nu}.$}

From now on, we fix $\nu \in \mathbb Z$ and focus on the checking the conditions 1 and 2 for the function $F^{nu}$ in Lemma \ref{L:TFAE}.
Since $\nu$ is fixed we will ommit the upper index $\nu$ in (\ref{E:Fnudecompdecomp})and denote $A_l(z_1)=A^{\nu}_l(z_1)$. Then we have
\begin{equation}\label{E:Fnu_bezindexa}
F^{\nu}(z_1,|z_2|)=\sum_{l=0}^{\infty}A_l(z_1)|z_2|^{2l}.
\end{equation}
We also know that
$$F^{\nu}(z_1,|z_2|)z_2^{\nu} \in C(\overline B^2).$$
\begin{lemma} \label{L:tfaetfae}
The conditions 1,2 of Lemma \ref{L:TFAE} are equivalent to the following: for any $r \in (0,1)$ the functions 
$A_l(\frac{z_1+c}{1+\overline c z_1})\frac{1}{(1+c \overline z_1)^{l}}$
extend continuosly  from $|z_1|=r$ to $|z_1|<r$ as a meromorphic function with the pole at $z=0$ of order at most $max(\nu,0).$
\end{lemma}
\pf
We have from (\ref{E:Fnu_bezindexa}):
$$(F^{\nu}\circ\omega_c)(z_1,z_2)=\sum\limits_{l=0}^{\infty} A_l(\frac{z_1+c}{1+\overline c z_1})\frac{|z_2|^{2l}}{|1+\overline c z_1|^{2l}}$$
Therefore for $|\lambda|=1$:
$$(F^{\nu}\circ\omega_c)(\lambda z)=\sum\limits_{l=0}^{\infty} A_l(\frac{\lambda z_1+c}{1+\overline c \lambda z_1})\frac{1}{(1+ c \overline {\lambda  z_1})^{l}}
\frac{1}{(1+ \overline c \lambda z_1)^{l}}|z_2|^{2l}.$$
Since $|z_2|$ is variable,  we conclude that the meromorphic extendibility of the function 
$$\lambda \mapsto (F^{\nu}\circ\omega_c)(\lambda z)$$ from 
the circle $|\lambda|=1$ to th edisc $|\lambda|<1$ with the pole  at $\lambda=0$ of order at most $ max(\nu,0).$
is equivalent to the meromorphic extendibility of the same type for the functions
$$\lambda \mapsto A_l(\frac{z_1+c_1}{1+\overline c_1 z_1})\frac{1}{(1+c \overline \lambda \overline z_1)^{l}}.$$
Here $|z_1|<1$ is arbitrary. Finally, the meromorphic extendibility in $\lambda$ from $|\lambda|=1$ can be obviously reworded as meromorphic extendibility in $z_1$ from the circles $|z_1|=r, \ 0<r<1. $
\begin{lemma}\label{L:A_lC1} For each $l$, the function 
\begin{equation}\label {E:B_l}
B_l(z_1):=A_l(z_1)(1-|z_1|^2)^{l + \nu}
\end{equation} 
develops in the open disc $|z_1|<1$ in power series in $z_1, \overline z_1$ and is continuous in the closed disc
$\{|z_1| \leq 1 \}.$ 
\end{lemma}
\pf The representation as a power series follows from the construction of the functions $B_l$ (\ref{E:B_l}), (\ref{E:Fnudecompdecomp}).
The function $F^{\nu}(z_1,|z_2|)$ is represented by a power series that converges in the ball $|z_1|^2+|z_1|^2 \leq 1.$ We can extend this representation in the complex space by replacing the variable $|z_2|^2$ by the complex variable $w=|z_2|^2.$
Then we obtain a function
\begin{equation}\label{E:FA}
F^{\nu}(z_1,w)=\sum_{l=0}^{\infty}A_l(z_1)w^l.
\end{equation}

For each fixed $z_1, \ |z_1|<1,$ the power series in $w$ converges in the disc $|w| \leq 1-|z_1|^2.$   The Taylor coefficient $A_l(z_1),$ corresponding to $w^l,$ expresses via the integral:
$$A_l(z_1)=\frac{l!}{2\pi i}\int\limits_{|\zeta|=1}\frac{F^{\nu}(z_1, (1-|z_1|^2)\zeta)}{(1-|z_1|^2)^{l}\zeta^{l+1}}d\zeta.$$
Then we have
\begin{equation}\label{E:Blz1}
B_l(z_1):=A_l(z_1)(1-|z_1|^2)^{l + \nu}=const \int\limits_{|\zeta|=1}F^{\nu}(z_1,(1-|z_1|^2)\zeta)(1-|z_1|^2)^{\nu}\frac{d\zeta}{\zeta^{l+1}}.
\end{equation}
The function $B_l$  is the sum of a power series and hence is continuous up to the circle $|z_1|=1$ because the function 
$$F^{\nu}(z_1,w)(1-|z_1|^2){\nu}$$ 
is continuous on the sphere $|z_1|^2+|w|^2=1$ by Lemma \ref{L:smooth}. 

\section{Meromorphic extension of functions $B_l^{\nu}(z_1)$ from hyperbolic circles}
We introduce the hyperbolic circles
$$H(c,r)=\{\frac{z_1+c}{1+\overline c z_1}:|z_1|=r.\}$$
The circle $H(c,r)$ is the image $H(c,r)=\omega_c(C_r)$ of the circle $C_r=\{|z_1|=r\}$ under the confimal automorihism $\omega_c(z)=(z+c)(1+\overline c z)^{-1})$ of the unit disc.
The point $c$ is the hyperbolic center of $H(c,r)$. The  Euclidean center and radius of the cirlce $H(c,r)$ are, correspondingly: 
$$e(c,r)=c\frac{1-r^2}{1-|c|^2r^2}, \ 
t(c,r)=r\frac{1-|c|^2}{1-|c|^2r^2}.$$
\begin{lemma}\label{L:B_l}
The  conditions 1,2 of Lemma \ref{L:TFAE} for the function $F^{\nu}$ and for the automorphism $\omega_c$ are equivalent to the following:

(*) for any $r \in (0,1), $ the functions $B_l,$ defined in Lemma \ref{L:A_lC1}, meromorphically extend from the hyperbolic circle $H(c,r)$ as a meromorphic function

with the only singular point- a pole at the Euclidean center $e(c,r)$ of order at most $max(\nu,0).$
\end{lemma}
\pf Due to Lemma \ref{L:tfaetfae}, it suffices to translate the conditions for  the functions $A_l$ formulated in Lemma \ref{L:tfaetfae} in terms of the functions 
$$ B_l(z_1):=A_l(z_1)(1-|z_1|^2)^{l + \nu}.$$
Consider the superposition of $B_l$ with the automorphism $\omega_c:$
$$B_l (\frac{z_1+c}{1+\overline c z_1})=\frac{(1-|z_1|^2)^{l+\nu}(1-|c|^2))^{l+ \nu}}{|1+\overline c z_1|^{2l + 2\nu}} A_l(\frac{z_1+c}{1+\overline c z_1}).$$ 
By Lemma \ref{L:tfaetfae}, the function 
$$M(z):=A_l(\frac{z_1+c}{1+\overline c z_1})\frac{1}{(1+c\overline z_1)^{l}}$$
extends from $|z_1|=r$ with the prescribed singularity at $z_1=0$ of order $\nu.$ 
We have on $|z_1|=r:$
$$B_l(\frac{z_1+c}{1+\overline c z_1})
=const \frac{1}{(1+\overline c z_1)^{l+\nu}}  M(z_1)\frac{1}{(1+ c\overline z_1)^{\nu}}.$$
On the circle $|z_1|=r$, we have 
$$\overline z_1=\frac{r^2}{z_1}$$ 
and therefore the meromorphic extension of $B_l$ into the disc $|z_1|<r$ is given by the function
$$
const \frac{1}{(1+\overline c z_1)^{l+\nu}}  M(z_1)\frac{z_1^{\nu}}{(z_1+ cr^2)^{\nu}}.
$$
Consider the case $\nu >0.$ The pole of $M(z_1)$ at $z_1=0$ is cancelled by the factor $z_1^{\nu}.$ Instead, the meromorphic extension develops a pole of the order $\leq \nu$ at the point
$$z_1=-cr^2.$$
It remains to notice that the automorphism $\omega_c$ sends the circle $|z_1|=1$ to the hyperbolic circle $H(c,r)$ and the point $-cr^2$ to the Euclidean center of $H(c,r)$: 
$$\omega_c(-cr^2)=\frac{-cr^2+c}{1 -|c|^2r^2}=e(c,r).$$ Therefore, the extension of $B_l\vert_{H(c,r)}$  develops a pole of order at most $\nu$ at the Euclidean center $e(c,r).$ The case $\nu <0$ is treated analogously, but in this case no poles appear.

\section{Reduction of Theorem \ref{T:mainmain} to characterization of polyanalytic functions in the unit disc}\label{S:reduction}

We will write $z$ instead of $z_1$ and $a,b$ instead of $a_2,b_1,$ correspondingly.
It will be convenient to introduce a temporary terminology. 

\begin{definition}\label{D:Hanu} Let $\nu \in \mathbb Z.$ We say that a function $E$ in the unit disc $\Delta$ satisfies the condition $(H,a,\nu)$ if
for any  hyperbolic circle $H(a,r), 0<r<1,$  the restriction
$f \vert_{H(a,r)}$ extends continuously in the disc bounded by $H(a,r)$ as a meromorphic function with the only singular point-a pole, of order at most $max(\nu,0),$ 
located at the Euclidean center $e(a,r)$ of the circle $H(a,r).$
\end{definition}

What we have proved in the previous sections can be summarized as follows:
\begin{itemize}
\item
The functions $B^{\nu}_l$ satisfy the conditions $(H,a,\nu)$ and $(H,b,\nu)$. 
\item Theorem \ref{T:mainmain}, and therefore the main result-Theorem \ref{T:main}, will be proved if we will be able to derive that $B_l^{\nu}$ satisfy
the condition $(H,c,\nu)$ for arbitrary point $c \in \Delta.$
\end{itemize}

It is easy to give sufficient conditions for functions to satisfy  the condition $(H,c,\nu)$ for all points $c \in \Delta.$ First of all, if $\nu \leq 0$ then any holomorphic function
$h$ in $\Delta$ is an example. If $\nu > 0$ then all polyanalytic functions of order $\nu$ satisfy the condition $(H,c,r)$:
\begin{lemma} \label{L:zbar}
Any function $E$ in $\Delta$, polyanalytic of order $\nu$ ,i.e. having the form
\begin{equation}\label{E:poly}
E(z)=h_0(z)+\overline zh_1(z)+\cdots \overline z^{\nu}h_{\nu}(z),
\end{equation}
where $h_j$ are holomorphic in $\Delta,$ satisfies the condition $(H,c,r)$ for all $c \in \Delta$ and $r \in (0,1).$
\end{lemma}
\pf
The set of circles $H(a,r)$ with arbitrary centers $c \in \Delta$ and radii $r \in (0,1)$  is just the set of all circles $C \subset \Delta.$
If $C=\{|z-e|=r\}$ then for $z \in C$ we have
$$\overline z=\overline e +\frac{r^2}{z-e}$$ and hence
\begin{equation}\label{E:polyanal}
E(z)=h_0(z)+ (  \overline e +\frac{r^2}{z-e})h_1(z)+\cdots + (\overline e +\frac{r^2}{z-e})^{\nu}h_{\nu}(z).
\end{equation}
Clearly, the meromorphic extension of $E(z)$ from $C$ in the disc $|z-e|<r$ develops a pole of order at most $\nu$  at the center $e.$ The real order of the pole depends on the order of zero
of the function $h_{\nu}(z)$ at $z=e.$ Of course, such cancellation can occur only for discrete set of the centers.

Observe also that the meromorphic extension $\overline z
=\overline e+ r^2/(z-e)$ inside the circle $|z-e|=r$ has no zero in the disc $|z-e|<r$ as long as $r>|a|,$ i.e. as long as the circle encloses the origin.

Taking into account Lemma \ref{L:zbar}, we have reduced the proof of Theorem \ref{T:mainmain} and hence, the proof of the main result Theorem \ref{T:main}, to the proof of the following:
\begin{proposition} \label{P:poly} The functions $B_l^{\nu}$ defined by (\ref{E:Fnudecompdecomp}), (\ref{E:B_l}) 
are polyanalytic of order $\leq \nu$ ($\nu>0)$ and analytic for $\nu \leq 0.$ Therefore $B_l^{\nu}$ satisfy the conditions $(H,c,\nu)$ for all $c \in \Delta$ and Theorem \ref{T:mainmain} and Theorem \ref{T:main} follow.
\end{proposition} 

In turn, Proposition \ref{P:poly} will follow from the characterization of polyanalytic functions which we are going to obtain now.
\section{Characterization of polyanalytic functions in the disc. Proof of Proposition \ref{P:poly}. End of the proofs of Theorem \ref{T:mainmain} and Theorem \ref{T:main}}

For simplicity, in this section, we will denote the complex variable in the unit disc by $z$ instead of $z_1$.  
\begin{theorem}\label{T:polyanalytic}
Let the power series $B(z)=\sum\limits_{n,m=0}^{\infty} b_{n,m}z^n\overline z^m$ converges in the unit disc $\{|z_1|<1\}$ and its sum $B(z_1)$ is continuous in the closed disc $\overline\Delta.$
Suppose that for a given $\nu \in \mathbb Z, $ the function $B$ satisfies the conditions $(H,a,\nu)$ and $(H,b,\nu)$ (Definition \ref{D:Hanu}) for some two points $a,b \in \Delta.$ Then
$B$ is polyanalytic in $\Delta$ of order at most $\nu$ (analytic in $\Delta$, if $\nu \leq 0$), i.e. has the form (\ref{E:polyanal}). 
\end{theorem}
Clearly, Proposition \ref{P:poly} is just Theorem \ref{T:polyanalytic} applied to the function $B=B^{\nu}_l.$
\begin{remark}
Theorem \ref{T:polyanalytic} was proved in \cite{APrep} for generic one-parameter families of circles. However, the proof in \cite{APrep}required stronger condition of regularity
for $B$ in th eclosed disc $\overline \Delta$ The  version presented here assumes pretty strong regularity (power decomposition) in the open  
disc but instead only continuity in the closed disc.  The proof is based on reduction to the case $\nu=0$ (holomorphic extension from  circles) 
and  referring to the result of \cite{Gl} which is obtained for this special families (concentric hyperbolic circles)  but instead under minimal regularity assumptions (just continuity). 
\end{remark}
\pf
Combining terms in the power series for $B(z)$ , we can represent the function $B(z)$ as
$$B(z)=\sum\limits_{k=0}^{\infty}h_k(z)\overline z^k,$$
where $h_k(z)$ are holomorphic functions in the unit disc $\Delta$, continuous in $\overline \Delta.$
Define
\begin{equation}\label{E:H}
H(z):=\frac{B(z)-h_0(z)-h_1(z)\overline z-\cdots-h_{\nu-1}(z)\overline z^{\nu-1}}{\overline z^{\nu}}=\sum\limits_{k=\nu}^{\infty}h_k(z)\overline z^{k-\nu}.
\end{equation}
The function $H$ is again real-analytic in $\Delta$ and continuous in $\overline \Delta.$

Now, we have seen in Lemma \ref{L:zbar} that the function $\overline z^{k}$ extends meromorphically from any circle with the pole of order $k$ at the Euclidean center, hence 
the difference $B(z)-h_0(z)-h_1(z)\overline z-\cdots-h_{\nu-1}\overline z^{\nu-1}$ extends from circles $H(a,r)$ and $H(b,r)$ with the poles at the Euclidean centers, of order at most $\nu.$

Let us analyze singularities of the meromorphic extensions of this difference after dividing by $\overline z^{\nu}.$
When the circle $H(a,r)$ encloses 0 then the meromorphic extension of $\overline z^{\nu}$ inside  $H(a,r)$ still has the pole of order $\nu$ at the Euclidean center $e(a,r)$ but has no zeros inside $H(a,r).$ 
This means that the poles of enumerator and denominator in (\ref{E:H}) cancel and no new poles appear.

Thus, if 0 is inside the circle $H(a,r)$, then the  function $H(z)$ extends holomorphically inside $H(a,r)$. This happens when the Euclidean radius and and center are related by
$t(a,r)>e(a,r).$ The holomorphic extendibilty means that the negative Fourier coefficients vanish:
$$
I(r):=\int\limits_{0}^{2\pi}H(e(a,r)+t(a,r)e^{i\phi})e^{im\phi}d\phi=0, \ m=1,\cdots,
$$
when $r \in (1-\varepsilon, 1]$ and $\varepsilon >0$ is sufficiently small.

However, $H(z)$ is real analytic and the functions $e(a,r), t(a,r)$ are real analytic for $r\in (0,1)$ and therefore $I(r)$  is real-analytic on $(0,1).$ Then by real-analyticity
$I(r)=0$ for all $r \in (0,1]$ which means 
that $H(z)$ analytically extends inside {\it all} circles $H(a,r)$ not only inside those enclosing $0.$ Analogously, $H$ extends holomorphically in all circles $H(b,r).$
Theorem 1.3 from \cite{Gl} asserts that then $H$ is holomorphic in $\Delta.$ 

Now, we have, by the construction of $H$:
$$B(z)=h_0(z)+h_1(z)\overline z +\cdots +h_{\nu-1}\overline z^{\nu-1}+H(z)\overline z^{\nu},$$
and now we know that $H(z)$ is holomorphic. Thus, $B$ is polyanalytic of order $\nu.$
Theorem \ref{T:polyanalytic} is proved.

According to the reduction summarized in Section \ref{S:reduction}, this proves Theorem \ref{T:mainmain} and hence proves Theorem \ref{T:main}.

\section{Proof of Theorem \ref{T:mainmainmain}}
The invariant Poisson integral $F=Pf$ coincides with $f$ on the  sphere $\partial B^2,$ hence, by the construction ( \ref{E:rewrite}) we have for $|z_1|^2+|z_2|^2=1:$
$$f(z_1,z_2)=F(z_1,z_2)=\sum\limits_{\nu=-\infty}^{\infty}F^{\nu}(z_1,|z_2|)z_2^{\nu}.$$

For $\nu <0,$ the function $F^{\nu}$ vanishes for $|z_1|=1, z_2=0$. 
Indeed, the coefficient functions $B^{\nu}_l$ extend from the circles $H(a,r),H(b,r)$ without singularities  hence they are 
holomorphic in the unit disc by Theorem \ref{T:polyanalytic} . Then $B^{\nu}_l=0$ since they vanish on the unit circle $|z_1|=1.$
Hence, $F^{\nu}=0$ for $\nu <0.$

Now let $\nu \geq 0$, Theorem \ref{T:polyanalytic} and  relations (\ref{E:Fnudecomp}), (\ref{E:B_l})
between the functions $F^{\nu}, A^{\nu}_l$ and $B^{\nu}_l$ yield for $z \in \partial B^2$:  
\begin{equation}\label{E:Fnu_last}
F^{\nu}(z)=\frac{1}{  (1-|z_1|^2)^{l+\nu} }  \sum\limits_{l=0}^{\infty}B^{\nu}_l(z_1) |z_2|^{2l}=\frac{1}{(1-|z_1|^2)^{\nu}}   \sum\limits_{l=0}^{\infty}B^{\nu}_l(z_1)  ,
\end{equation}
where 
\begin{equation}\label{E:Bnu}
B^{\nu}_l(z_1)=
h^{\nu}_{l,0}(z_1)+... + h^{\nu}_{l,\nu}(z_1)\overline z_1^{\nu}.
\end{equation}
Write 
$$\overline z_1=\frac{(|z_1|^2-1)+1}{z}.$$
Substitute this expression for $\overline z_1$ in the (\ref{E:Bnu}). We obtain a new representation of the form
\begin{equation}\label{E:B_last}
B^{\nu}_l(z_1)=g^{\nu}_{l,0}(z_1)+ g^{\nu}_{l,1}(z_1)(1-|z_1|^2)+...+ g^{\nu}_{l,\nu}(z_1)(1-|z_1|^2)^{\nu},
\end{equation}
where $g^l_j(z)$ are analytic in $\Delta$ but maybe with poles at $z=0.$ 

However, singularities at $z=0$ are impossible because then the meromorphic extensions of $B^{\nu}_l(z_1)$ in any of the circles $H(a,r), H(b,r)$ would have poles at $z=0$ 
which is not the case if $e(a,r) \neq 0, \ e(b,r) \neq 0.$
Therefore, the singularities of the functions $g_j(z_1)$  at $z=0$ are removable and they are holomorphic in $\Delta.$

Substitute (\ref{E:B_last}) into (\ref{E:Fnu_last}):
$$F^{\nu}(z_1,|z_2|)=\sum\limits_{l=0}^{\infty}\frac{g^{\nu}_{l,0}(z_1)}{(1-|z_1|^2)^{\nu}}+\frac{g^{\nu}_{l,1}(z_1)}{(1-|z_1|^2)^{\nu-1}}+...+g^{\nu}_{l,\nu}(z_1).$$
Summation in the index $l$ implies the representation
$$F^{\nu}(z_1,|z_2|)=\frac{h^\nu_\nu(z_1)}{(1-|z_1|^2)^{\nu}}+\frac{h^{\nu}_{\nu-1}(z_1)}{(1-|z_1|^2)^{\nu-1}}+...+h^{\nu}_{0}(z_1),$$
where $h^{\nu}_j(z_1)$ are holomorphic.
The corresponding series of $g^{\nu}_{l,j}$ in $l$  converge because $(1-|z_1|^2){\nu} F^{\nu}$ is the sum of a power series in $B^2$ and the new representation is nothing but regrouping terms. 
Then the representation \ref{E:ff} for the function $f$ on $\partial B^2$ follows if one replaces $1-|z_1|^2$ by $|z_2|^2,$ using the equation of the sphere. It remains to observe that the continuity of each 
term $F^{\nu}(z_1,|z_2|)z_2^{\nu}$ in the power series for the Poisson integral $F=Pf$
requires that $h^{\nu}_j(z_1)=0, z_1 \in \overline\Delta,$ whenever $2j \geq  \nu$ and hence the summation is performed  only as long as $2j <\nu.$ 
Theorem \ref{T:mainmainmain} is proved.

\section*{Acknowledgements}
This work was partially supported by the grant from ISF (Israel Science Foundatrion) 688/08.
Some of this research was done as a part of European Science Foundation Networking Program HCAA.

\end{document}